\documentclass{amsart}

\usepackage{latexsym, amsfonts, amsmath, amssymb, amscd, graphicx}

\newtheorem{theorem}{Theorem}[section]

\newtheorem{corollary}[theorem]{Corollary}
\newtheorem{prop}[theorem]{Proposition}

\theoremstyle{definition}

\newtheorem{remark}[theorem]{Remark}
\newtheorem{df}[theorem]{Definition}

\newcommand{\diag}{\mathrm{diag}}
\newcommand{\sg}{\mathrm{Sym}}
\newcommand{\Sl}{\mathfrak{sl}}
\newcommand{\Psl}{\mathfrak{psl}}
\newcommand{\Pgl}{\mathfrak{pgl}}

\newcommand{\So}{\mathfrak{so}}
\newcommand{\Sp}{\mathfrak{sp}}

\newcommand{\GL}{\mathrm{GL}}
\newcommand{\PGL}{\mathrm{PGL}}
\newcommand{\Ort}{\mathrm{O}}
\newcommand{\SP}{\mathrm{Sp}}
\newcommand{\A}{\mathcal{A}}
\newcommand{\AI}{\mathcal{A}^\mathrm{(I)}}
\newcommand{\AII}[1]{\mathcal{A}^{\mathrm{(II_{#1})}}}
\newcommand{\B}{\mathcal{B}}
\newcommand{\C}{\mathcal{C}}
\newcommand{\D}{\mathcal{D}}
\newcommand{\M}{\mathcal{M}}

\newcommand{\bg}{{\overline{g}}}
\newcommand{\bt}{{\overline{t}}}
\newcommand{\be}{{\overline{e}}}
\newcommand{\bG}{{\overline{G}}}
\newcommand{\bT}{{\overline{T}}}
\newcommand{\ZZ}{\mathbb{Z}}

\newcommand{\FF}{\mathbb{F}}

\newcommand{\OO}{\mathbb{O}}

\newcommand{\sks}{\mathcal{K}}
\newcommand{\wh}[1]{\widehat{#1}}
\newcommand{\wt}[1]{\widetilde{#1}}

\newcommand{\vphi}{\varphi}
\newcommand{\vp}{\varphi}
\newcommand{\veps}{\varepsilon}
\newcommand{\chr}[1]{\mathrm{char}\,#1}
\newcommand{\ad}[1]{\mathrm{ad}\,#1}

\newcommand{\Aut}{\mathrm{Aut}\,}
\newcommand{\Lie}[1]{{#1}^{(-)}}
\newcommand{\Antaut}{\overline{\mathrm{Aut}}\,}
\newcommand{\AAut}{\mathbf{Aut}\,}
\newcommand{\Der}{\mathrm{Der}\,}

\newcommand{\im}{\mathrm{im}\,}
\newcommand{\lspan}[1]{\mathrm{Span}\,\{#1\}}
\newcommand{\matr}[1]{\left(\begin{smallmatrix}#1\end{smallmatrix}\right)}
\newcommand{\bee}[1]{\begin{equation}\label{#1}}
\newcommand{\ene}{\end{equation}}
\newcommand{\ot}{\otimes}
\newcommand{\sgn}{\mathrm{sgn}}
\newcommand{\ld}{\ldots}

\begin{document}

\title[Classification of Gradings on Lie algebras]{Classification of Group Gradings\\ on Simple Lie algebras of Types A, B, C and D}

\author[Bahturin]{Yuri Bahturin}
\address{Department of Mathematics and Statistics\\Memorial
University of Newfoundland\\ St. John's, NL, A1C5S7, Canada}
\email{bahturin@mun.ca}

\author[Kochetov]{Mikhail Kochetov}
\address{Department of Mathematics and Statistics\\Memorial
University of Newfoundland\\ St. John's, NL, A1C5S7, Canada}
\email{mikhail@mun.ca}

\thanks{The first author acknowledges support by by NSERC grant \# 227060-04. The second author acknowledges support by NSERC Discovery Grant \# 341792-07.}

\begin{abstract}
For a given abelian group $G$, we classify the isomorphism classes of $G$-gradings on the simple Lie algebras of types $\A_n$ ($n\geq 1$), $\B_n$ ($n\geq 2$), $\C_n$ ($n\geq 3$) and $\D_n$ ($n>4$), in terms of numerical and group-theoretical invariants. The ground field is assumed to be algebraically closed of characteristic different from $2$. 
\end{abstract}

\maketitle


\section{Introduction}\label{intro}

Let $U$ be an algebra (not necessarily associative) over a field $\FF$ and let $G$ be an abelian group, written multiplicatively.

\begin{df}\label{group_grad}
A {\em $G$-grading} on $U$ is a vector space decomposition
\[
U=\bigoplus_{g\in G} U_g
\]
such that 
\[
U_g U_h\subset U_{gh}\;\mbox{ for all }\;g,h\in G.
\]
$U_g$ is called the {\em homogeneous component} of degree $g$. The {\em support} of the $G$-grading is the set
\[
\{g\in G\;|\;U_g\neq 0\}.
\] 
\end{df}

\begin{df}\label{iso_grad}
We say that two $G$-gradings, $U=\bigoplus_{g\in G} U_g$ and $U=\bigoplus_{g\in G} U'_g$, are {\em isomorphic} if there exists an algebra automorphism $\psi:U\to U$ such that
\[
\psi(U_g)=U'_{g}\;\mbox{ for all }\;g\in G,
\]
i.e., $U=\bigoplus_{g\in G} U_g$ and $U=\bigoplus_{g\in G} U'_g$ are isomorphic as $G$-graded algebras.
\end{df}

The purpose of this paper is to classify, for a given abelian group $G$, the isomorphism classes of $G$-gradings on the classical simple Lie algebras of types $\A_n$ ($n\geq 1$), $\B_n$ ($n\geq 2$), $\C_n$ ($n\geq 3$) and $\D_n$ ($n>4$), in terms of numerical and group-theoretical invariants. Descriptions of such gradings were obtained in \cite{BShZ,BZA,antaut,BKM,BK}, but the question of distinguishing non-isomorphic gradings was not addressed in those papers. Also, A. Elduque \cite{Eld09} has recently found a counterexample to \cite[Proposition 6.4]{BZA}, which was used in the description of gradings on Lie algebras of type $\A$. The fine gradings (i.e., those that cannot be refined) on Lie algebras of types $\A$, $\B$, $\C$ and $\D$ (including $\D_4$) have been classified, up to equivalence, in \cite{Eld09} over algebraically closed fields of characteristic zero. For a discussion of the difference between classification up to equivalence and classification up to isomorphism see \cite{Ksur}. The two kinds of classification cannot be easily obtained from each other.

We will assume throughout this paper that the ground field $\FF$ is algebraically closed. We will usually assume that $\chr{\FF}\neq 2$ and in one case also $\chr{\FF}\neq 3$. We obtain a description of gradings in type $\A$ without using \cite[Proposition 6.4]{BZA} and with methods simpler that those in \cite{BKM,BK}. We also obtain invariants that allow us to distinguish among non-isomorphic gradings in types $\A$, $\B$, $\C$ and $\D$.

The paper is structured as follows. In Section \ref{section_matr} we recall the description of $G$-gradings on a matrix algebra $R=M_n(\FF)$ and determine when two such gradings are isomorphic (Theorem \ref{matr_main}). We also obtain a canonical form for an anti-automorphism of $R$ that preserves the grading and restricts to an involution on the identity component $R_e$ (Theorem \ref{form_antaut}). In particular, this allows us to classify (up to isomorphism) the pairs $(R,\vphi)$ where $R=M_n(\FF)$ is $G$-graded and $\vphi$ is an involution that preserves the grading (Corollary \ref{matr_with_inv}). In Section \ref{section_transfer} we use affine group schemes to show how one can reduce the classification of $G$-gradings on classical simple Lie algebras to the classification of $G$-gradings on $R=M_n(\FF)$ and of the pairs $(R,\vphi)$ where $\vphi$ is an involution or an anti-automorphism satisfying certain properties. In Section \ref{section_A} we obtain a classification of $G$-gradings on simple Lie algebras of type $\A$ --- see Theorem \ref{A_main}. Finally, in Section \ref{section_BCD} we state a classification of $G$-gradings on simple Lie algebras of types $\B$, $\C$ and $\D$ (except $\D_4$) --- see Theorem \ref{BCD_main}, which is an immediate consequence of Corollary \ref{matr_with_inv}.


\section{Gradings on Matrix Algebras}\label{section_matr}

Let $R=M_n(\FF)$ where $\FF$ is an algebraically closed field of arbitrary characteristic. Let $G$ be an abelian group. A description of $G$-gradings on $R$ was obtained in \cite{BSZ,BZnc,surgrad}. In this section we restate that description in a slightly different form and obtain invariants that allow us to distinguish among non-isomorphic gradings. Criteria for isomorphism of the so-called ``elementary'' gradings (see below) on matrix algebras $M_n(\FF)$ and on the algebra of finitary matrices were obtained in \cite{Cha} and \cite{BZfin}, respectively.

We start with gradings $R=\bigoplus_{g\in G}R_g$ with the property $\dim R_g\leq 1$ for all $g\in G$. As shown in the proof of \cite[Theorem 5]{BSZ}, $R$ is then a graded division algebra, i.e., any nonzero homogeneous element is invertible in $R$. Consequently, the support $T\subset G$ of the grading is a subgroup. Following  \cite{Eld09}, we will call such $R=\bigoplus_{g\in G}R_g$ a {\em division grading} (the terms used in \cite{BSZ,BZnc,surgrad} and in \cite{HPP} are ``fine gradings'' and ``Pauli gradings'', respectively). Note that since $R\cong \FF^\sigma T$ is semisimple, $\chr{\FF}$ does not divide $n^2=|T|$.

For each $t\in T$, let $X_t$ be a nonzero element in the component $R_t$. Then 
\[
X_u X_v=\sigma(u,v)X_{uv}
\] 
for some nonzero scalar $\sigma(u,v)$. Clearly, the function $\sigma:T\times T\to\FF^{\times}$ is a $2$-cocycle, and the $G$-graded algebra $R$ is isomorphic to the twisted group algebra $\FF^\sigma T$ (with its natural $T$-grading regarded as a $G$-grading). Rescaling the elements $X_t$ corresponds to replacing $\sigma$ with a cohomologous cocycle. Let 
\[
\beta_\sigma(u,v):=\frac{\sigma(u,v)}{\sigma(v,u)}.
\]
Then $\beta=\beta_\sigma$ depends only on the class of $\sigma$ in $H^2(T,\FF^\times)$ and $\beta:T\times T\to \FF^\times$ is an {\em alternating bicharacter}, i.e., it is multiplicative in each variable and has the property $\beta(t,t)=1$ for all $t\in T$. 

Clearly, $X_uX_v=\beta(u,v)X_vX_u$. Since the centre $Z(R)$ is spanned by the identity element, $\beta$ is {\em nondegenerate} in the sense that $\beta(u,t)=1$ for all $u\in T$ implies $t=e$. Conversely, if $\sigma$ is a $2$-cocycle such that $\beta_\sigma$ is nondegenerate, then $\FF^\sigma T$ is a semisimple associative algebra whose centre is spanned by the identity element, so $\FF^\sigma T$ is isomorphic to $R$. Therefore, the isomorphism classes of division $G$-gradings on $R=M_n(\FF)$ with support $T\subset G$ are in one-to-one correspondence with the classes $[\sigma]\in H^2(T,\FF^\times)$ such that $\beta_\sigma$ is nondegenerate.

The classes $[\sigma]$ and the corresponding gradings on $R$ can be found explicitly as follows. As shown in the proof of \cite[Theorem 5]{BSZ}, there exists a decomposition of $T$ into the direct product of cyclic subgroups:
\begin{equation}\label{decomp_T}
T=H_1'\times H_1''\times\cdots\times H_r'\times H_r''
\end{equation} 
such that $H_i'\times H_i''$ and $H_j'\times H_j''$ are $\beta$-orthogonal for $i\neq j$, and $H_i'$ and $H_i''$ are in duality by $\beta$. Denote by $\ell_i$ the order of $H_i'$ and $H_i''$. If we pick generators $a_i$ and $b_i$ for $H_i'$ and $H_i''$, respectively, then $\veps_i:=\beta(a_i,b_i)$ is a primitive $\ell_i$-th root of unity, and all other values of $\beta$ on the elements $a_1,b_1,\ldots,a_r,b_r$ are 1. Pick elements $X_{a_i}\in R_{a_i}$ and $X_{b_i}\in R_{b_i}$ such that $X_{a_i}^{\ell_i}=X_{b_i}^{\ell_i}=1$. Then we obtain an isomorphism $\FF^\sigma T\to M_{\ell_1}(\FF)\ot\cdots\ot M_{\ell_r}(\FF)$ defined by
\begin{equation}
X_{a_i}\mapsto I\ot\cdots I\ot X_i\ot I\ot\cdots I\quad\mbox{and}\quad X_{b_i}\mapsto I\ot\cdots I\ot Y_i\ot I\ot\cdots I,
\end{equation}
where
\begin{equation}\label{Pauli}
X_i=\begin{bmatrix} 
\veps_i^{n-1} & 0                   & 0           & \ldots      & 0       & 0\\ 
0             & \veps_i^{n-2}       & 0           & \ldots      & 0       & 0\\ 
\ldots        &                     &             &             &         &  \\[3pt]
0             & 0                   & 0           & \ldots      & \veps_i & 0\\ 
0             & 0                   & 0           & \ldots      & 0       & 1
\end{bmatrix}\mbox{ and }
Y_i=\begin{bmatrix} 
0 & 1 & 0 & \ldots & 0 & 0\\ 
0 & 0 & 1 & \ldots & 0 & 0\\ 
\ldots & & & & \\[3pt] 
0 & 0 & 0 & \ldots & 0 & 1\\
1 & 0 & 0 & \ldots & 0 & 0
\end{bmatrix}
\end{equation}
are in the $i$-th factor, $M_{\ell_i}(\FF)$.

It follows that the class $[\sigma]\in H^2(T,\FF^\times)$, and hence the isomorphism class of the $G$-graded algebra $\FF^\sigma T$, is uniquely determined by $\beta=\beta_\sigma$. Conversely, since the relation $X_uX_v=\beta(u,v)X_vX_u$ does not change when we rescale $X_u$ and $X_v$, the values of $\beta$ are determined by the $G$-grading. We summarize our discussion in the following

\begin{prop}\label{class_division}
There exist division $G$-gradings on $R=M_n(\FF)$ with support $T\subset G$ if and only if $\chr{\FF}$ does not divide $n$ and $T\cong\ZZ_{\ell_1}^2\times\cdots\times\ZZ_{\ell_r}^2$ where $\ell_1\cdots\ell_r=n$. The isomorphism classes of division $G$-gradings with support $T$ are in one-to-one correspondence with nondegenerate alternating bicharacters $\beta:T\times T\to\FF^\times$.\hfill{$\square$}
\end{prop}

We also note that taking
\[
X_{(a_1^{i_1},b_1^{j_1},\ldots,a_r^{i_r},b_r^{j_r})}=X_{a_1}^{i_1}X_{b_1}^{j_1}\cdots X_{a_r}^{i_r}X_{b_r}^{j_r},
\]
we obtain a representative of the cohomology class $[\sigma]$ that is multiplicative in each variable, i.e., it is a bicharacter (not alternating unless $T$ is the trivial subgroup). In what follows, we will always assume that $\sigma$ is chosen in this way.

\begin{df}\label{division_model}
A concrete representative of the isomorphism class of division $G$-graded algebras with support $T$ and bicharacter $\beta$ can be obtained as follows. First decompose $T$ as in (\ref{decomp_T}) and pick generators $a_1,b_1,\ldots,a_r,b_r$. Then define a grading on $M_{\ell_i}(\FF)$ by declaring that $X_i$ has degree $a_i$ and $Y_i$ has degree $b_i$, where $X_i$ and $Y_i$ are given by (\ref{Pauli}) and $\veps_i=\beta(a_i,b_i)$. Then $M_{\ell_1}(\FF)\ot\cdots\ot M_{\ell_r}(\FF)$ with tensor product grading is a representative of the desired class. We will call any representative obtained in this way a {\em standard realization}.
\end{df}

If $R$ has a division grading, then its structure is quite rigid. Any automorphism of the graded algebra $R$ must send $X_t$ to a scalar multiple of itself, hence it is given by $X_t\mapsto\lambda(t)X_t$ where $\lambda:T\to\FF^\times$ is a character of $T$. Since $\beta$ is nondegenerate, it establishes an isomorphism between $T$ and $\wh{T}$. It follows that the automorphism of $R$ corresponding to $\lambda$ is given by $X\mapsto X_t^{-1}XX_t$ where $t\in T$ is determined by $\beta(u,t)=\lambda(u)$ for all $u\in T$. 

It follows from \cite[Lemma 6.1]{BZA} that the graded algebra $R$ admits anti-auto\-mor\-phisms only when $T$ is an elementary $2$-group (and hence $\chr{\FF}\neq 2$). In this case, we can regard $T$ as a vector space over the field of order 2 and think of $\sigma(u,v)$ as a bilinear form on $T$. Hence $\sigma(t,t)$ is a quadratic form, and $\beta(u,v)$ is the polar bilinear form for $\sigma(t,t)$. Note that $\sigma(t,t)$ depends on the choice of $\sigma$, so it is not an invariant of the graded algebra $R$. In fact, any quadratic form with polar form $\beta(u,v)$ can be achieved by changing generators $a_i,b_i$ in the $i$-th copy of $\ZZ_2^2$. However, once we fix a standard realization of $R$, $\sigma(t,t)$ is uniquely determined. Following the usual convention regarding quadratic forms, we will denote $\sigma(t,t)$ by $\beta(t)$ so that $\beta(u,v)=\beta(uv)\beta(u)\beta(v)$. Note that 
\begin{equation}\label{inv_division}
X^\beta=\beta(u)X\quad\mbox{for all}\quad X\in R_u, u\in T
\end{equation} 
is an involution of the graded algebra $R$. Hence any anti-automorphism of $R$ is given by $X\mapsto X_t^{-1}X^\beta X_t$ for a suitable $t\in T$. In the standard realization of $R$ as $M_2(\FF)^{\ot r}$, the involution $\beta$ is given by matrix transpose. We summarize the above discussion for future reference:

\begin{prop}\label{antaut_division}
Suppose $R=M_n(\FF)$ has a division $G$-grading with support $T\subset G$ and bicharacter $\beta$. Then the mapping that sends $t\in T$ to the inner automorphism $X\mapsto X_t^{-1}XX_t$ is an isomorphism between $T$ and the group of automorphisms $\Aut_G(R)$ of the graded algebra $R$. The graded algebra $R$ admits anti-automorphisms if and only if $T$ is an elementary $2$-group. If this is the case, then, in any standard realization of $R$, the mapping $X\mapsto {}^tX$ is an involution of the graded algebra $R$. This involution can be written in the form (\ref{inv_division}), where $\beta:T\to\{\pm 1\}$ is a quadratic form. The bicharacter $\beta(u,v)$ is the polar bilinear form associated to $\beta$. The group $\Antaut_G(R)$ of automorphisms and anti-automorphisms of the graded algebra $R$ is equal to $\Aut_G(R)\times \langle\beta\rangle$. In particular, any anti-automorphism of the graded algebra $R$ is an involution, given by $X\mapsto X_t^{-1}X^\beta X_t$ for a uniquely determined $t\in T$.\hfill{$\square$}
\end{prop}

We now turn to general $G$-gradings on $R$. As shown in \cite{BSZ,BZnc,surgrad}, there exist graded unital subalgebras $C$ and $D$ in $R$ such that $D\cong M_\ell(\FF)$ has a division grading, $C\cong M_k(\FF)$ has an {\em elementary grading} given by a $k$-tuple $(g_1,\ldots,g_k)$ of elements of $G$:
\[
C_g=\lspan{E_{ij}\;|\;g_i^{-1}g_j=g}\quad\mbox{for all}\quad g\in G,
\]
where $E_{ij}$ is a basis of matrix units in $C$, and we have an isomorphism $C\ot D\to R$ given by $c\ot d\mapsto cd$. Moreover, the intersection of the support $\{g_i^{-1}g_j\}$ of the grading on $C$ and the support $T$ of the grading on $D$ is equal to $\{e\}$. 

Without loss of generality, we may assume that the $k$-tuple has the form 
\[
(g_1^{(k_1)},\ldots,g_s^{(k_s)})
\] 
where the elements $g_1,\ldots,g_s$ are pairwise distinct and we write $g^{(q)}$ for $\underbrace{g,\ldots,g}_{q\;{\rm times}}$. 

It is important to note that the subalgebras $C$ and $D$ are not uniquely determined. We are now going to obtain invariants of the graded algebra $R$. The partition $k=k_1+\cdots+k_s$ gives a block decomposition of $C$. Let $e_i$ be the block-diagonal matrix $\diag(0,\ldots,I_{k_i},\ldots,0)$ where $I_{k_i}$ is in the $i$-th position, $i=1,\ldots,s$. Consider the Peirce decomposition of $C$ corresponding to the orthogonal idempotents $e_1,\ldots,e_s$: $C_{ij}=e_iCe_j$. We will write $C_i$ instead of $C_{ii}$ for brevity. Then the identity component is 
\[
R_e=C_1\ot I\oplus\cdots\oplus C_s\ot I.
\]
It follows that the idempotents $e_1,\ldots,e_s$ and the (non-unital) subalgebras $C_1,\ldots,C_s$ of $R$ are uniquely determined (up to permutation). It is easy to verify that the centralizer of $R_e$ in $R$ is equal to $e_1\ot D\oplus\cdots\oplus e_s\ot D$. Hence the (non-unital) subalgebras $D_i:=e_i\ot D$ of $R$ are uniquely determined (up to permutation). All $D_i$ are isomorphic to $D$ as $G$-graded algebras, so the isomorphism class of $D$ is uniquely determined. This gives us invariants $T$ and $\beta$ according to Proposition \ref{class_division}. However, there is no canonical way to choose the isomorphisms of $D$ with $D_i$. According to Proposition \ref{antaut_division}, the possible choices are parameterized by $t_i\in T$, $i=1,\ldots,s$. If we fix isomorphisms $\eta_i:D\to D_i$, then each Peirce component $R_{ij}=e_iRe_j$ becomes a $D$-bimodule by setting $d\cdot r=\eta_i(d)r$ and $r\cdot d=r\eta_j(d)$ for all $d\in D$ and $r\in R_{ij}$. Taking $\eta_i(d)=e_i\ot d$ for all $d\in D$, we recover the subspaces $C_{ij}$ for $i\neq j$ as the centres of these bimodules:
\[
C_{ij}=\{r\in R_{ij}\;|\;d\cdot r=r\cdot d\quad\mbox{for all}\quad d\in D\}.
\]
Also, the subalgebra $D$ of $R$ can be identified:
\[
D=\{\eta_1(d)+\cdots+\eta_s(d)\;|\;d\in D\}.
\]
If we replace $\eta_i$ by $\eta'_i(d)=\eta_i(X_{t_i}^{-1}dX_{t_i})$, then we get $C'_{ij}=\eta_i(X_{t_i}^{-1})C_{ij}\eta_j(X_{t_j})$. 
Let $C'=C_1\oplus\cdots\oplus C_s\oplus\bigoplus_{i\neq j}C'_{ij}$ and $D'=\{\eta'_1(d)+\cdots+\eta'_s(d)\;|\;d\in D\}$. Then $C'$ and $D'$ are graded unital subalgebras of $R$. Let $\Psi=e_1\ot X_{t_1}+\cdots+e_s\ot X_{t_s}$. Then $\Psi$ is an invertible matrix and the mapping $\psi(X)=\Psi^{-1}X\Psi$ is an automorphism of the (ungraded) algebra $R$ that sends $C$ to $C'$ and $D$ to $D'$. The restriction of $\psi$ to $D$ preserves the grading, whereas the restriction of $\psi$ to $C$ sends homogeneous elements of degree $g_i^{-1}g_j$ to homogeneous elements of degree $t_i^{-1}g_i^{-1}g_jt_j$ (i.e., ``shifts'' the grading in the $(i,j)$-th Peirce components by $t_i^{-1}t_j$). We conclude that the $G$-grading of $R$ associated to the $k$-tuple $(g_1^{(k_1)},\ldots,g_s^{(k_s)})$ is isomorphic to the $G$-grading associated to the $k$-tuple $((g_1t_1)^{(k_1)},\ldots,(g_st_s)^{(k_s)})$. Finally, we note that the cosets $g_i^{-1}g_jT$ are uniquely determined by the $G$-graded algebra $R$, because they are the supports of the grading on the Peirce components $R_{ij}$ ($i\neq j$). We have obtained an irredundant classification of $G$-gradings on $R$. 

To state the result precisely, we introduce some notation. Let 
\[
\kappa=(k_1,\ldots,k_s)\quad\mbox{where}\; k_i\;\mbox{are positive integers.}
\]
We will write $|\kappa|$ for $k_1+\cdots+k_s$ and $e_i$, $i=1,\ldots,s$, for the orthogonal idempotents in $M_{|\kappa|}(\FF)$ associated to the block decomposition determined by $\kappa$. Let
\[
\gamma=(g_1,\ldots,g_s)\quad\mbox{where}\; g_i\in G\;\mbox{are such that}\; g_i^{-1}g_j\notin T\;\mbox{for all}\;i\neq j.
\] 

\begin{df}\label{equiv_1}
We will write $(\kappa,\gamma)\sim (\wt{\kappa},\wt{\gamma})$ if $\kappa$ and $\wt{\kappa}$ have the same number of components $s$ and there exist an element $g\in G$ and a permutation $\pi$ of the symbols $\{1,\ldots,s\}$ such that $\wt{k}_i=k_{\pi(i)}$ and $\wt{g}_i\equiv g_{\pi(i)}g\pmod{T}$, for all $i=1,\ldots,s$.
\end{df}

\begin{df}\label{datum_1}
Let $D$ be a standard realization of division $G$-graded algebra with support $T\subset G$ and bicharacter $\beta$. Let $\kappa$ and $\gamma$ be as above. Let $C=M_{|\kappa|}(\FF)$. We endow the algebra $M_{|\kappa|}(D)=C\ot D$ with a $G$-grading by declaring the degree of $U\ot d$ to be $g_i^{-1}tg_j$ for all $U\in e_iCe_j$ and $d\in D_t$. We will denote this $G$-graded algebra by $\M(G,T,D,\kappa,\gamma)$. By abuse of notation, we will also write $\M(G,T,\beta,\kappa,\gamma)$,
since the isomorphism class of $D$ is uniquely determined by $\beta$.
\end{df}

\begin{theorem}\label{matr_main}
Let $\FF$ be an algebraically closed field of arbitrary characteristic. Let $G$ be an abelian group. Let $R=\bigoplus_{g\in G}R_g$ be a grading of the matrix algebra $R=M_n(\FF)$. Then the $G$-graded algebra $R$ is isomorphic to some $\M(G,T,\beta,\kappa,\gamma)$ where $T\subset G$ is a subgroup, $\beta:T\times T\to\FF^\times$ is a nondegenerate alternating bicharacter, $\kappa$ and $\gamma$ are as above with $|\kappa|\sqrt{|T|}=n$. Two $G$-graded algebras $\M(G,T_1,\beta_1,\kappa_1,\gamma_1)$ and $\M(G,T_2,\beta_2,\kappa_2,\gamma_2)$ are isomorphic if and only if $T_1=T_2$, $\beta_1=\beta_2$ and $(\kappa_1,\gamma_1)\sim(\kappa_2,\gamma_2)$.\hfill{$\square$}
\end{theorem}

\begin{remark}\label{monomial_iso}
In fact, it follows from the above discussion that, for any permutation $\pi$ as in Definition \ref{equiv_1}, there exists an isomorphism from $\M(G,T,\beta,\wt{\kappa},\wt{\gamma})$ to $\M(G,T,\beta,\kappa,\gamma)$ that sends $\wt{e}_i$ to $e_{\pi(i)}$. We can construct such an isomorphism explicitly in the following way. Let $P=P_\pi$ be the block matrix with $I_{k_i}$ in the $(i,\pi(i))$-th positions and $0$ elsewhere (i.e., the block-permutation matrix corresponding to $\pi$). Pick $t_i\in T$ such that $\wt{g}_i=g_{\pi(i)} t_{\pi(i)} g$ and let $B$ be the block-diagonal matrix $e_1\ot X_{t_1}+\cdots+e_s\ot X_{t_s}$. Then the map $X\mapsto (BP)X(BP)^{-1}$ has the desired properties. We will refer to isomorphisms of this type as {\em monomial}.
\end{remark} 

Let $\sg(s)$ be the group of permutations on $\{1,\ldots,s\}$. Let $\Aut(\kappa,\gamma)$ be the subgroup of $\sg(s)$ that consists of all $\pi$ such that, for some $g\in G$, we have $k_i=k_{\pi(i)}$ and $g_i\equiv g_{\pi(i)}g\pmod{T}$ for all $i=1,\ldots,s$.  


\begin{prop}\label{aut_grading}
The group of automorphisms $\Aut_G(R)$ of the graded algebra $R=\M(G,T,\beta,\kappa,\gamma)$ is an extension of $\Aut(\kappa,\gamma)$ by $\PGL_\kappa(\FF)\times\Aut_G(D)$ where
\[
\PGL_\kappa(\FF)=\big(\GL_{\kappa_1}(\FF)\times\cdots\times\GL_{\kappa_s}(\FF)\big)/\FF^\times,
\]
where $\FF^\times$ is identified with nonzero scalar matrices.
\end{prop}

\begin{proof}
Any $\psi\in\Aut_G(R)$ leaves the identity component $R_{e}$ invariant and hence permutes the idempotents $e_1,\ldots,e_s$. This gives a homomorphism $f:\Aut_G(R)\to\sg(s)$. Looking at the supports of the Peirce components, we see that $f(\psi)\in\Aut(\kappa,\gamma)$. Conversely, any element of $\Aut(\kappa,\gamma)$ is in $\im f$ by Remark \ref{monomial_iso}, since it comes from a monomial automorphism of the graded algebra $R$. Finally, any $\psi\in\ker f$ leaves $C_i$ and $D_i$ invariant and hence is given by $\psi(X)=\Psi^{-1}X\Psi$ where $\Psi=B_1\ot Q_1\oplus\cdots\oplus B_s\ot Q_s$ for some $B_i\in\GL_{\kappa_i}(\FF)$ and $Q_i\in D$. In view of Proposition \ref{antaut_division}, we may assume that $Q_i=X_{t_i}$ for some $t_i\in T$. It is easy to see that $\psi$ preserves the grading if and only if $t_1=\ldots=t_s$. The result follows.
\end{proof}

In order to classify gradings on Lie algebras of types $\B$, $\C$ and $\D$, we will need to study involutions on $G$-graded matrix algebras. A description of such involutions was given in \cite{antaut}. Here we will slightly simplify that description and obtain invariants that will allow us to distinguish among isomorphism classes. We start with a more general situation, which we will need for the classification of gradings in type $\A$.

\begin{df}
Let $G$ be an abelian group and let $U=\bigoplus_{g\in G} U_g$ be a $G$-graded algebra. We will say that an anti-automorphism $\vphi$ of $U$ is {\em compatible} with the grading if $\vphi(U_g)=U_g$ for all $g\in G$. If $U_1$ and $U_2$ are $G$-graded algebras and $\vphi_1$ and $\vphi_2$ are anti-automorphisms on $U_1$ and $U_2$, respectively, compatible with the grading, then we will say that $(U,\vphi_1)$ and $(U,\vphi_2)$ are {\em isomorphic} if there exists an isomorphism $\psi:U_1\to U_2$ of $G$-graded algebras such that $\vphi_1=\psi^{-1}\vphi_2\psi$.
\end{df}

Suppose $R=M_n(\FF)$ is $G$-graded and there exists an anti-automorphism $\vphi$ compatible with the grading and such that $\vphi^2|_{R_e}=id$. Then $\vphi$ leaves some of the components of $R_e$ invariant and swaps the remaining components in pairs. Without loss of generality, we may assume that $e_i$ are $\vphi$-invariant for $i=1,\ldots, m$ and not $\vphi$-invariant for $i>m$. It will be convenient to change the notation and write $e_{m+1}',e_{m+1}'',\ldots,e_k',e_k''$ so that $\vphi$ swaps $e_i'$ and $e_i''$ for $i>m$. (Thus the total number of orthogonal idempotents in question is $2k-m$.) It will also be convenient to distinguish $\vphi$-invariant idempotents of even and odd rank. Thus we assume that $e_1,\ldots,e_\ell$ have odd rank and $e_{\ell+1},\ldots,e_m$ have even rank. We will change the notation for $\kappa$ and $\gamma$ accordingly:
\begin{equation}\label{datum_2_kappa}
\kappa=(q_1,\ldots,q_\ell,2q_{\ell+1},\ldots,2q_m,q_{m+1},q_{m+1}\ldots,q_k,q_k)
\end{equation}
where $q_i$ are positive integers with $q_1,\ldots,q_\ell$ odd, and
\begin{equation}\label{datum_2_gamma}
\gamma=(g_1,\ldots,g_\ell,g_{\ell+1},\ldots,g_m, g_{m+1}',g_{m+1}'',\ldots,g_k',g_k'')
\end{equation}
where $g_i\in G$ are such that $g_i^{-1}g_j\notin T$ for all $i\neq j$.

As shown in \cite{BZA,antaut}, the existence of the anti-automorphism $\vphi$ places strong restrictions on the $G$-grading. First of all, note that the centralizer of $R_e$ in $R$, which is equal to $D_1\oplus\cdots\oplus D_m\oplus D_{m+1}'\oplus D_{m+1}''\oplus\cdots\oplus D_k'\oplus D_k''$, is $\vphi$-invariant. Since $e_1,\ldots,e_m$ are $\vphi$-invariant and belong to $D_1,\ldots,D_m$, respectively, we see that $D_1,\ldots, D_m$ are also $\vphi$-invariant. By a similar argument, $\vphi$ swaps $D_i'$ and $D_i''$ for $i>m$. Each of the $D_i$, $D_i'$ and $D_i''$ is an isomorphic copy of $D$, so we see that $D$ admits an anti-automorpism. By Proposition \ref{antaut_division}, $T$ must be an elementary $2$-group and we have a standard realization $D\cong M_2(\FF)^{\ot r}$. 

Since $\vphi$ preserves the $G$-grading and $\vphi(e_iRe_j)=e_jRe_i$ for $i,j\leq m$, the supports of these two Peirce components must be equal, which gives $g_i^{-1}g_j\equiv g_j^{-1}g_i\pmod{T}$ for $i,j\leq m$. Similarly, $\vphi(e_i'Re_j'')=e_j'Re_i''$ implies $(g_i')^{-1}g_j''\equiv (g_j')^{-1}g_i''\pmod{T}$ for $i,j>m$. Also, $\vphi(e_iRe_j')=e_j''Re_i$ implies $g_i^{-1}g_j'\equiv (g_j'')^{-1}g_i\pmod{T}$ for $i\leq m$ and $j>m$. These conditions can be summarized as follows:
\begin{equation}\label{compat_datum_2_gamma_mod}
g_1^2\equiv\ldots\equiv g_m^2\equiv g_{m+1}'g_{m+1}''\equiv\ldots\equiv g_k'g_k''\pmod{T}.
\end{equation}
If $\gamma$ satisfies (\ref{compat_datum_2_gamma_mod}), then we have
\[
g_1^2t_1=\ldots=g_m^2t_m=g_{m+1}'g_{m+1}''t_{m+1}=\ldots=g_k'g_k''t_k
\]
for some $t_1,\ldots,t_k\in T$. We can replace the $G$-grading by an isomorphic one so that $\gamma$ satisfies
\begin{equation}\label{compat_datum_2_gamma}
g_1^2t_1=\ldots=g_m^2t_m=g_{m+1}'g_{m+1}''=\ldots=g_k'g_k''.
\end{equation}
Indeed, it suffices to replace $g_i''$ by $g_i''t_i$, $i=m+1,\ldots,k$ (which does not change the cosets mod $T$).

\begin{theorem}\label{form_antaut}
Let $\FF$ be an algebraically closed field, $\chr{\FF}\neq 2$. Let $G$ be an abelian group. Let $R=\M(G,T,\beta,\kappa,\gamma)$. Assume that $R$ admits an anti-automorphism $\vphi$ that is compatible with the grading and satisfies $\vphi^2|_{R_e}=id$. Write $\kappa$ and $\gamma$ in the form (\ref{datum_2_kappa}) and (\ref{datum_2_gamma}), respectively. Then $T$ is an elementary $2$-group and $\gamma$ satisfies (\ref{compat_datum_2_gamma_mod}). Up to an isomorphism of the pair $(R,\vphi)$, $\gamma$ satisfies (\ref{compat_datum_2_gamma}) for some $t_1,\ldots,t_m\in T$ and $\vphi$ is given by $\vphi(X)=\Phi^{-1}({}^tX)\Phi$ for all $X\in R$, where matrix $\Phi$ has the following block-diagonal form:
\begin{equation}\label{formula_antaut}
\Phi=\sum_{i=1}^\ell{I_{q_i}\ot X_{t_i}}\oplus \sum_{i=\ell+1}^m{S_i\ot X_{t_i}}\oplus
\sum_{i=m+1}^k\begin{pmatrix}0&I_{q_i}\\\mu_i I_{q_i}&0\end{pmatrix}\ot I
\end{equation}
where, for $i=\ell+1,\ldots,m$, each $S_i$ is either $I_{2q_i}$ or $\matr{0&I_{q_i}\\-I_{q_i}&0}$. 
\end{theorem}

\begin{proof}
There exists an invertible matrix $\Phi$ such that $\vphi$ is given by $\vphi(X)=\Phi^{-1}({}^tX)\Phi$ for all $X\in R$. Recall that conjugating $\vphi$ by the automorphism $\psi(X)=\Psi^{-1}X\Psi$ replaces matrix $\Phi$ by ${}^t\Psi\Phi\Psi$, i.e., $\Phi$ is transformed as the matrix of a bilinear form. 

Recall that we fixed the idempotents 
\begin{equation}\label{idempotents}
e_1,\ld,e_{\ell},e_{\ell+1},\ld,e_m,e_{m+1}^{\prime},e_{m+1}^{\prime\prime},\ld,e_{k}^{\prime},e_{k}^{\prime\prime}.
\end{equation}
It is also convenient to introduce $e_i=e_i'+e_i''$ for $i=m+1,\ldots, k$.

Following the proof of \cite[Lemma 6 and Proposition 1]{antaut}, we see that, up to an automorphism of the $G$-graded algebra $R$, $\Phi$ has the following block-diagonal form --- in agreement with the idempotents given by (\ref{idempotents}):
\[
\Phi=\sum_{i=1}^\ell{S_iY_i\ot Q_i\oplus \sum_{i=\ell+1}^mS_iY_i\ot Q_i}\oplus\sum_{i=m+1}^kS_iY_i\ot Q_i.
\]
(This is formula (20) of just cited paper, rewritten according to our present notation.) For $i=1,\ld,m$, the matrix $Y_i$ is in the centralizer of the simple algebra $C_i$, i.e., has the form $Y_i=\xi_iI_{q_i}$. For $i=m+1,\ld,k$, the matrix $Y_i$ is in the centralizer of the semisimple algebra $C_i^{\prime}\oplus C_i^{\prime\prime}$, i.e., has the form  $Y_i=\diag(\eta_i I_{q_i},\xi_i I_{q_i})$. Each $Q_i$ is in $D_i$, and the map $X\mapsto Q_i^{-1}({}^tX) Q_i$ is an anti-automorphism of $D$. Hence, by Proposition \ref{antaut_division}, each $Q_i$ is, up to a scalar multiple, of the form $X_{t_i}$, for an appropriate choice of $t_i\in T$. The scalar can be absorbed in $Y_i$. Finally,  the matrix $S_i$ is $I_{q_i}$ for $i=1,\ld,\ell$, either $I_{2q_i}$ or $\matr{0&I_{q_i}\\-I_{q_i}&0}$ for $i=\ell+1,\ldots,m$, and $\matr{0&I_{q_i}\\I_{q_i}&0}$ for $i=m+1,\ldots,k$. This allows us to rewrite the above formula as follows:
\[
\Phi=\sum_{i=1}^\ell\xi_iI_{q_i}\ot X_{t_i}\oplus \sum_{i=\ell+1}^m\xi_iS_i\ot X_{t_i}\oplus
\sum_{i=m+1}^k\begin{pmatrix}0&\xi_i I_{q_i}\\\eta_iI_{q_i}&0\end{pmatrix}\ot X_{t_i}.
\]
Here $\xi_i,\eta_i$ are some nonzero scalars. If we now apply the inner automorphism of the graded algebra $R$ given by the matrix 
$P=\frac{1}{\sqrt{\xi_1}}e_1\ot I+\cdots+\frac{1}{\sqrt{\xi_k}}e_k\ot I$, then $\vphi$ is transformed to the anti-automorphism given by the following matrix (which we again denote by $\Phi$):
\[
\Phi=\sum_{i=1}^\ell I_{q_i}\ot X_{t_i}\oplus \sum_{i=\ell+1}^mS_i\ot X_{t_i}\oplus
\sum_{i=m+1}^k\begin{pmatrix}0& I_{q_i}\\\mu_iI_{q_i}&0\end{pmatrix}\ot X_{t_i},
\]
for an appropriate set of nonzero scalars $\mu_{m+1},\ld,\mu_k$. It can be easily verified (and is shown in the proof of \cite[Theorem 3]{antaut}) that $t_1,\ldots,t_k$ satisfy the following condition:
$g_1^2t_1=\ldots=g_m^2t_m=g_{m+1}'g_{m+1}''t_{m+1}=\ldots=g_s'g_s''t_s$.

Finally, the inner automorphism $\psi(X)=\Psi^{-1}X\Psi$ of $R$ where 
\[
\Psi^{-1}=e_1\ot I+\cdots+e_m\ot I+e_{m+1}'\ot I+e_{m+1}''\ot X_{t_{m+1}}+\cdots+e_k'\ot I+e_k''\ot X_{t_k}
\] 
sends the $G$-grading to the one given by 
\[
(g_1,\ldots,g_m,g_{m+1}',g_{m+1}''t_{m+1},\ldots,g_k',g_k''t_k)
\]
and transforms $\vphi$ to the anti-automorphism given by a matrix of form (\ref{formula_antaut}). 
\end{proof}

If $\vphi$ is an involution on $R$, then one can get rid of the parameters $\mu_{m+1},\ldots,\mu_s$, and the selection of $S_{\ell+1},\ldots,S_m$ is uniquely determined. Set $\sgn(\vphi)=1$ if $\vphi$ is a transpose involution and $\sgn(\vphi)=-1$ if $\vphi$ is a symplectic involution on $R$. Similarly, set $\sgn(S_i)=1$ if ${}^tS_i=S_i$ and $\sgn(S_i)=-1$ if ${}^tS_i=-S_i$. We restate the main result of \cite{antaut} in our notation (and setting $t_{m+1}=\ldots=t_k=e$):

\begin{theorem}\cite[Theorem 3]{antaut}\label{form_inv}
Under the conditions of Theorem \ref{form_antaut}, assume that $\vphi^2=id$. Then, up to an isomorphism of the pair $(R,\vphi)$, $\gamma$ satisfies (\ref{compat_datum_2_gamma}) for some $t_1,\ldots,t_m\in T$ and $\vphi$ is given by  $\vphi(X)=\Phi^{-1}({}^tX)\Phi$ for all $X\in R$, where matrix $\Phi$ has the following block-diagonal form:
\begin{equation}\label{formula_inv}
\Phi=\sum_{i=1}^\ell{I_{q_i}\ot X_{t_i}}\oplus \sum_{i=\ell+1}^m{S_i\ot X_{t_i}}\oplus\sum_{i=m+1}^k S_i\ot I
\end{equation}
where
\begin{itemize} 
\item for $i=\ell+1,\ldots,m$, each $S_i$ is either $I_{2q_i}$ or $\matr{0&I_{q_i}\\-I_{q_i}&0}$, and
\item for $i=m+1,\ldots,s$, all $S_i$ are either $\matr{0&I_{q_i}\\I_{q_i}&0}$ or $\matr{0&I_{q_i}\\-I_{q_i}&0}$ 
\end{itemize}
such that the following condition is satisfied:
\begin{equation}\label{compat_inv}
\begin{split}
\sgn(\vphi)&=\beta(t_1)=\ldots=\beta(t_\ell)\\
&=\beta(t_{\ell+1})\sgn(S_{\ell+1})=\ldots=\beta(t_m)\sgn(S_m)\\
&=\sgn(S_{m+1})=\ldots=\sgn(S_k).
\end{split}
\end{equation}
Conversely, if $\gamma$ satisfies (\ref{compat_datum_2_gamma}) and condition (\ref{compat_inv}) holds, then $\Phi$ defines an involution of the type indicated by $\sgn(\vphi)$ on the $G$-graded algebra $R$.\hfill{$\square$}
\end{theorem}

It is convenient to introduce the following notation (for $m>0$):
\begin{equation}\label{datum_2_tau}
\tau=(t_1,\ldots,t_m).
\end{equation}

Note that for the elements $t_1,\ldots,t_m$ in (\ref{compat_datum_2_gamma}), the ratios $t_i^{-1}t_j$ are uniquely determined by the cosets of $g_1,\ldots,g_m$ mod $T$, so it is sufficient to specify only one $t_i$ to find $\tau$.

\begin{df}
We will say that $\gamma$ is {\em $*$-admissible} if it satisfies (\ref{compat_datum_2_gamma_mod}) and, for some $t_1,\ldots,t_\ell\in T$, we have $g_1^2t_1=\ldots=g_\ell^2t_\ell$ and 
\begin{equation}\label{inv_admissible}
\beta(t_1)=\ldots=\beta(t_\ell).
\end{equation}
(If $\ell\leq 1$, then condition (\ref{inv_admissible}) is automatically satisfied.) 
\end{df}

\begin{df}\label{datum_2}
Let $T\subset G$ be an elementary $2$-group (of even rank) with a nondegenerate alternating bicharacter $\beta$. Suppose $\gamma$ is $*$-admissible, and $\gamma$ and $\tau$ satisfy (\ref{compat_datum_2_gamma}) and (\ref{inv_admissible}). If $\ell>0$, let $\delta$ be the common value of $\beta(t_1),\ldots,\beta(t_\ell)$. If $\ell=0$, select $\delta\in\{\pm 1\}$ arbitrarily. Consider $R=\M(G,T,\beta,\kappa,\gamma)$. Let $\Phi$ be the matrix given by (\ref{formula_inv}) where the matrices $S_i$ are selected so that equation  (\ref{compat_inv}) holds with $\sgn(\vphi)=\delta$. Then, by Theorem \ref{form_inv}, $\vphi(X)=\Phi^{-1}({}^tX)\Phi$ is an involution on $R$ that is compatible with the grading.
We will denote $(R,\vphi)$ defined in this way by $\M^*(G,T,\beta,\kappa,\gamma,\tau,\delta)$. (Here $\tau$ is empty if $m=0$.)
\end{df}

\begin{df}\label{equiv_2}
Referring to Definition \ref{datum_2}, we will write $(\kappa,\gamma,\tau)\approx (\wt{\kappa},\wt{\gamma},\wt{\tau})$ if $\kappa$ and $\wt{\kappa}$ have the same number of components of each type, i.e., the same values of $\ell$, $m$ and $k$, and there exist an element $g\in G$ and a permutation $\pi$ of the symbols $\{1,\ldots,k\}$ preserving the sets $\{1,\ldots,\ell\}$, $\{\ell+1,\ldots,m\}$ and $\{m+1,\ldots,k\}$ such that $\wt{q}_i=q_{\pi(i)}$ for all $i$, $\wt{g}_i\equiv g_{\pi(i)}g\pmod{T}$ for all $i=1,\ldots,m$, $\{\wt{g}'_i,\wt{g}''_i\}\equiv\{g'_{\pi(i)}g,g''_{\pi(i)}g\}\pmod{T}$ for all $i=m+1,\ldots,k$, and
\begin{itemize}
\item if $m>0$, then $\wt{t}_i=t_{\pi(i)}$ for all $i=1,\ldots,m$;
\item if $m=0$, then $\wt{g}'_i\wt{g}''_i=g'_{\pi(i)}g''_{\pi(i)}g^2$ for some (and hence all) $i=1,\ldots,k$.
\end{itemize}
In the case $m=0$, $\tau$ is empty, so we may write $(\kappa,\gamma)\approx (\wt{\kappa},\wt{\gamma})$.
\end{df}

\begin{corollary}\label{matr_with_inv}
Let $\chr{\FF}\neq 2$ and $R=\M(G,T,\beta,\kappa,\gamma)$. Then the $G$-graded algebra $R$ admits an involution if and only if $T$ is an elementary $2$-group and $\gamma$ is $*$-admissible. If $\vphi$ is an involution on $R$, then $(R,\vphi)$ is isomorphic to some $\M^*(G,T,\beta,\kappa,\gamma,\tau,\delta)$ where $\delta=\sgn(\vphi)$. Two $G$-graded algebras with involution, $\M^*(G,T_1,\beta_1,\kappa_1,\gamma_1,\tau_1,\delta_1)$ and $\M^*(G,T_2,\beta_2,\kappa_2,\gamma_2,\tau_2,\delta_2)$, are isomorphic if and only if $T_1=T_2$, $\beta_1=\beta_2$, $(\kappa_1,\gamma_1,\tau_1)\approx(\kappa_2,\gamma_2,\tau_2)$ and $\delta_1=\delta_2$.
\end{corollary}

\begin{proof}
The first two statements are a combination of Theorems \ref{form_antaut} and \ref{form_inv}. It remains to prove the last statement. 

Let $R_1=\M(G,T_1,\beta_1,\kappa_1,\gamma_1)$,  $R_2=\M(G,T_2,\beta_2,\kappa_2,\gamma_2)$ and let $\vphi_1$ and $\vphi_2$ be the corresponding involutions. Suppose $T_1=T_2$, $\beta_1=\beta_2$, and $(\kappa_1,\gamma_1,\tau_1)\approx(\kappa_2,\gamma_2,\tau_2)$.  Then, by Theorem \ref{matr_main}, there exists an isomorphism of $G$-graded algebras  $\psi:R_1\to R_2$. By Remark \ref{monomial_iso}, $\psi$ can be chosen to be a monomial isomorphism associated to the permutation $\pi$ in Definition \ref{equiv_2}. The matrix of the involution $\psi^{-1}\vphi_2\psi$ on $R_1$ is then obtained from the matrix of $\vphi_2$ by permuting the blocks on the diagonal so that they align with the corresponding blocks of $\vphi_1$ and possibly multiplying some of the blocks by $-1$ (the extra condition for the case $m=0$ in Definition \ref{equiv_2} guarantees that the second tensor factor in each block remains $I$). If $\delta_1=\delta_2$, then  $\psi^{-1}\vphi_2\psi$ can be transformed to $\vphi_1$ by an automorphism of the $G$-graded algebra $R_1$ (see the proof of Theorem \ref{form_antaut}). 

Conversely, suppose there exists an isomorphism $\psi:(R_1,\vphi_1)\to(R_2,\vphi_2)$. First of all, $\delta_1$ and $\delta_2$ are determined by the type of involution (transpose or symplectic), so $\delta_1=\delta_2$. By Theorem \ref{matr_main}, we also have $T_1=T_2$, $\beta_1=\beta_2$, $(\kappa_1,\gamma_1)\sim(\kappa_2,\gamma_2)$. The partitions of $\kappa_1$ and $\kappa_2$ according to $\{1,\ldots,\ell\}$, $\{\ell+1,\ldots,m\}$ and $\{m+1,\ldots,k\}$ are determined by $\vphi_1$ and $\vphi_2$, hence they must correspond under $\psi$. At the same time, for some $g\in G$, the cosets of $\gamma_1gT$ and $\gamma_2T$ must correspond under $\psi$ up to switching $g'_i$ with $g''_i$ ($i>m$). In the case $m>0$, by Proposition \ref{antaut_division}, $\tau_1$ and $\tau_2$ are uniquely determined by the restrictions of $\vphi_1$ and $\vphi_2$ to $D_1,\ldots,D_m$ and hence must match under the permutation determined by $\psi$. Therefore, in this case $(\kappa_1,\gamma_1,\tau_1)\approx(\kappa_2,\gamma_2,\tau_2)$. It remains to consider the case $m=0$. Looking at the description of the automorphism group given by Proposition \ref{aut_grading}, we see that $\psi=\psi_0\alpha$ where $\psi_0$ is a monomial isomorphism and $\alpha$ is in $\PGL_{\kappa_1}(\FF)\times\Aut_G(D)$.  The action of $\psi_0$ on $\vphi_2$ leads to the permitation of blocks and the replacement of the second tensor factor $I$ by $X_{t_0}$ for some $t_0\in T$. Then $\alpha$ must transform $\psi_0^{-1}\vphi_2\psi_0$ to $\vphi_1$. The effect of $\alpha$ on one block is the following (we omit subscripts to simplify notation):
\[\begin{split}
&\left(\begin{pmatrix}{}^tA&0\\0&{}^tB\end{pmatrix}\ot{}^tX_{u}\right)\left(\begin{pmatrix}0&I\\ \veps I&0\end{pmatrix}\ot X_{t_0}\right)\left(\begin{pmatrix}A&0\\0&B\end{pmatrix}\ot X_u\right)\\
&=\pm\begin{pmatrix}0&{}^tAB\\ \veps{}\,^tBA&0\end{pmatrix}\ot X_{t_0}.
\end{split}\]
We see that $\alpha$ cannot change $t_0$. It follows that $t_0=e$ and $(\kappa_1,\gamma_1)\approx(\kappa_2,\gamma_2)$.
\end{proof}


\section{Correspondence between Lie Gradings and Associative Gradings}\label{section_transfer}

Let $U$ be an algebra and let $G$ be a group. Then a $G$-grading on $U$ is equivalent to a structure of an $\FF G$-comodule algebra (see e.g. \cite{Mont} for background). If we assume that $U$ is finite-dimensional and $G$ is abelian and finitely generated, then the comodule structure can be regarded as a morphism of (affine) algebraic group schemes $G^D\to\AAut(U)$ where $G^D$ is the Cartier dual of $G$ and $\AAut(G)$ is the automorphism group scheme of $U$ (see e.g. \cite{Wh} for background). Two $G$-gradings are isomorphic if and only if the corresponding morphisms $G^D\to\AAut(U)$ are conjugate by an automorphism of $U$. Note also that, if $U$ is finite-dimensional, then we may always assume without loss of generality that $G$ is finitely generated (just replace $G$ by the subgroup generated by the support of the grading).

If $\chr{\FF}=0$, then $G^D=\wh{G}$, the algebraic group of characters on $G$, and $\AAut(G)=\Aut(G)$, the algebraic group of automorphisms. If $\chr{\FF}=p>0$, then we can write $G=G_0\times G_1$ where $G_0$ has no $p$-torsion and $G_1$ is a $p$-group. Hence $G^D=\wh{G_0}\times G_1^D$, where $\wh{G_0}$ is smooth and $G_1^D$ is finite and connected. The algebraic group $\wh{G_0}$ (which is equal to $\wh{G}$) acts on $U$ as follows: 
\[
\chi\ast X=\chi(g)X\quad\mbox{for all}\;X\in U_g\;\mbox{and}\;g\in G. 
\]

The group scheme $\AAut(U)$ contains the group $\Aut(U)$ as the largest smooth subgroupscheme. The tangent Lie algebra of $\AAut(U)$ is $\Der(U)$, so $\AAut(U)$ is smooth if and only if $\Der(U)$ equals the tangent Lie algebra of the group $\Aut(U)$.

We will be interested in the following algebras: $M_n(\FF)$, $\Psl_n(\FF)$, $\So_n(\FF)$ and $\Sp_n(\FF)$. In all these cases the automorphism group scheme is smooth, i.e., coincides with the algebraic group of automorphisms (regarded as a group scheme). Indeed, for the associative algebra $R=M_n(\FF)$, it is well-known that $\Aut(R)=\PGL_n(\FF)$ and $\Der(R)=\Pgl_n(\FF)$. For the Lie algebra $L=\So_n(\FF)$ ($n\geq 5$, $n\neq 8$) or $\Sp_n(\FF)$ ($n\geq 4$), it is known that every automorphism of $L$ is the conjugation by an element of $\Ort_n(\FF)$ or $\SP_n(\FF)$, respectively --- see \cite{Jac} for the case $\chr{\FF}=0$ and \cite{Sel} for the case $\chr{\FF}=p$ ($p\neq 2$). In particular, every automorphism of $L$ is the restriction of an automorphism of $R$. Similarly, every derivation of $L$ is the restriction of a derivation of $R$ (see e.g. \cite{BBCM3}). 

Let $\vphi$ be the involution of $R$ such that $L=\sks(R,\vphi)$, the space of skew-symmetric elements with repsect to $\vphi$. Then the projectivizations of the groups $\Ort_n(\FF)$ and $\SP_n(\FF)$ are equal to $\Aut(R,\vphi)$, and their tangent algebras are equal to $\Der(R,\vphi)$. Hence the restriction map $\theta:\Aut(R,\vphi)\to\Aut(L)$ is a surjective homomorphism of algebraic groups such that $d\theta:\Der(R,\vphi)\to\Der(L)$ is also surjective. It follows that $\AAut(L)$ is smooth. Since $L$ generates $R$ as an associative algebra, both $\theta$ and $d\theta$ are also injective. Hence  $\theta:\Aut(R,\vphi)\to\Aut(L)$ is an isomorphism of algebraic groups. For $G$-gradings this means the following. Clearly, if $R=\bigoplus_{g\in G}R_g$ is a grading that is compatible with $\vphi$, then the restriction $L_g=R_g\cap L$ is a grading of $L$. Since $\theta:\Aut(R,\vphi)\to\Aut(L)$ is an isomorpism and the automorphism groups are equal to the automorphism group schemes, the restriction map gives a bijection between the isomorphism classes of $G$-gradings on $L$ and the $\Aut(R,\vphi)$-orbits on the set of $\vphi$-compatible $G$-gradings on $R$. The orbits correspond to isomorphism classes of pairs $(R,\vphi)$ where $R=M_n(\FF)$ is $G$-graded and $\vphi$ is an involution on $R$ that is compatible with the grading.

The case of $L=\Psl_n(\FF)$ is more complicated. We have a homomorphism of algebraic groups $\theta:\Aut(R)\to\Aut(L)$ given by restriction and passing to cosets modulo the centre. It is well-known that this homomorphism is not surjective for $n\geq 3$, because the map $X\mapsto-\,{}^tX$ is not an automorphism of the associative algebra $R$, but it is an automorphism of the Lie algebra $\Lie{R}$ and hence induces an automorphism of $L$. Let $\Antaut(R)$ be the group of automorphisms and anti-automorphisms of $R$. Then we can extend $\theta$ to a homomorphism $\Antaut(R)\to\Aut(L)$ by sending an anti-automorphism $\vphi$ of $R$ to the map induced on $L$ by $-\vphi$. This extended $\theta$ is surjective for any $n\geq 3$ if $\chr{\FF}\neq 2,3$ (see \cite{Sel}) and for any $n>3$ if $\chr{\FF}=3$ (see \cite{BBCM3}). It is easy to verify that $\theta$ and $d\theta$ are injective and hence $\theta$ is an isomorphism of algebraic groups (see e.g. \cite[Lemma 5.3]{BK}). It is shown in \cite{BBCM3} that, under the same assumptions on $\chr{\FF}$, every derivation of $L$ is induced by a derivation of $R$. It follows that $\AAut(L)$ is smooth, i.e., $\AAut(L)=\Aut(L)$.

Now let $L=\bigoplus_{g\in G}L_g$ be a $G$-grading and let $\alpha:G^D\to\Aut(L)$ be the corresponding morphism. Then we have a morphism $\wt{\alpha}:=\theta^{-1}\alpha:G^D\to\Antaut(R)$, which gives a $G$-grading $R=\bigoplus_{g\in G}R_g$ on the Lie algebra $\Lie{R}$. The two gradings are related in the following way: $L_g=(R_g\cap[R,R])$ mod $Z(R)$.

Set $\Lambda=\wt{\alpha}^{-1}(\Aut(R))$. Then $\Lambda$ is a subgroupscheme of $G^D$ of index at most $2$. Moreover, since $G_1^D$ is connected, it is mapped by $\wt{\alpha}$ to $\Aut(R)$ and hence is contained in $\Lambda$. We have two possibilities: either $\Lambda=G^D$ or $\Lambda$ has index $2$. Following \cite{BZA}, we will say that the $G$-grading on $L$ has {\em Type I} in the first case and has {\em Type II} in the second case. In Type I, the $G$-grading corresponding to $\wt{\alpha}$ is a grading of $R$ as an associative algebra. In Type II, we consider $\Lambda^\perp$, which is a subgroup of order $2$ in $G$. Let $h$ be the generator of this subgroup.
\begin{remark}
For the readers more familiar with the language of Hopf algebras, there is an alternative way to define the element $h$. The Hopf algebra $\FF[\Antaut(R)]$ of regular functions on the algebraic group $\Antaut(R)$ has a group-like element $f$ defined by $f(\psi)=1$ if $\psi$ is an automorphism and $f(\psi)=-1$ if $\psi$ is an anti-automorphism. The morphism of group schemes $\wt{\alpha}:G^D\to\Antaut(R)$ corresponds to a homomorphism of Hopf algebras $\FF[\Antaut(R)]\to\FF G$. The element $h$ is the image of $f$ under this homomorphism.
\end{remark}
Let $\bG=G/\langle h\rangle$. Then the restriction $\wt{\alpha}:\Lambda\to\Aut(R)$ corresponds to the coarsening of the $G$-grading on $R$ given by the quotient map $G\to\bG$:
\[
R=\bigoplus_{\bg\in\bG}R_\bg\quad\mbox{where}\;R_\bg=R_g\oplus R_{gh}.
\] 
This $\bG$-grading is a grading of $R$ as an associative algebra. The $G$-grading on $\Lie{R}$ can be recovered as follows. Fix $\chi\in\wh{G_0}=\wh{G}$ such that $\chi(h)=-1$. Then $\chi$ acts on $R$ as $-\vphi$ where $\vphi$ is an anti-automorphism preserving the $\bG$-grading. Then we have
\[
R_g=\{X\in R_\bg\;|\;-\vp(X)=\chi(g)X\}=\{-\vp(X)+\chi(g)X\;|\;X\in R_\bg\}.
\]
Thus we obtain 1) a bijection between the isomorphism classes of $G$-gradings on $L$ of Type I and the $\Antaut(R)$-orbits on the set of $G$-gradings on $R$ and 2) a bijection between the isomorphism classes of $G$-gradings on $L$ of Type II and $\Antaut(R)$-orbits on the set of pairs $(R,\vp)$ where $R=M_n(\FF)$ is $\bG$-graded and $\vphi$ is an anti-automorphism on $R$ that is compatible with the $\bG$-grading and has the property $\vp^2(X)=\chi^2\ast X$ for all $X\in R$.

\begin{remark}
If $n=2$, then $\theta:\Aut(R)\to\Aut(L)$ is an isomorphism, so there are no gradings of Type II.
\end{remark}


\section{Gradings on Lie Algebras of Type $\A$}\label{section_A}

Let $L=\Psl_n(\FF)$ and $R=M_n(\FF)$, where $\chr{\FF}\neq 2$ and, for $n=3$, also $\chr{\FF}\neq 3$. Let $L=\bigoplus_{g\in G}L_g$ be a grading of $L$ by an abelian group $G$. As discussed in the previous section, this grading belongs to one of two types. Gradings of Type I are induced from $G$-grading on the associative algebra $R$, which have been classified in Theorem \ref{matr_main}. 

\begin{df}\label{type_I}
Let $R=\M(G,T,\beta,\kappa,\gamma)$ and let $L_g=(R_g\cap[R,R])$ mod $Z(R)$. We will denote the $G$-graded algebra $L$ obtained in this way as $\AI(G,T,\beta,\kappa,\gamma)$.
\end{df}

Now assume that we have a grading of Type II. Then there is a distinguished element $h\in G$ of order 2. 
Let $\bG=G/\langle h\rangle$. Then the $G$-grading on $L$ is induced from a $G$-grading on the Lie algebra $\Lie{R}$ that is obtained by refining a $\bG$-grading 
$R=\bigoplus_{\bg\in\bG}R_\bg$ on the associative algebra $R$. Let $R=\M(\bG,\bT,\beta,\kappa,\gamma)$ as a $\bG$-graded algebra. The refinement is obtained using the action of any character $\chi\in\wh{G}$ with $\chi(h)=-1$, and the result does not depend on the choice of $\chi$. So we fix $\chi\in\wh{G}$ such that $\chi(h)=-1$. 

Set $\vp(X)=-\chi\ast X$ for all $X\in R$. Then $\vp$ is an anti-automorphism of the $\bG$-graded algebra $R$. Moreover, $\vp^2(X)=\chi^2\ast X$. Since $\chi^2(h)=1$, we can regard $\chi^2$ as a character on $\bG$ and hence its action on $X\in R_{\bg}$ is given by $\chi^2\ast X = \chi^2(\bg)X$. In particular, $\vp^2|_{R_\be}=id$. By Theorem \ref{form_antaut}, $\bT$ is an elementary $2$-group, $\kappa$ is given by (\ref{datum_2_kappa}) and $\gamma$ is given by (\ref{datum_2_gamma}) with bars over the $g$'s. We may also assume that $\gamma$ satisfies 
\begin{equation}\label{compat_A_gamma}
\bg_1^2\bt_1=\ldots=\bg_m^2\bt_m=\bg_{m+1}'\bg_{m+1}''=\ldots=\bg_s'\bg_s''
\end{equation}
for some $\bt_1,\ldots,\bt_m\in\bT$, and $\vp$ is given by $\vp(X)=\Phi^{-1}({}^tX)\Phi$ where
\begin{equation}\label{formula_antaut_restated}
\Phi=\sum_{i=1}^\ell{I_{q_i}\ot X_{\bt_i}}\oplus \sum_{i=\ell+1}^m{S_i\ot X_{\bt_i}}\oplus
\sum_{i=m+1}^k\begin{pmatrix}0&I_{q_i}\\\mu_i I_{q_i}&0\end{pmatrix}\ot I,
\end{equation}
where $\mu_i$ are nonzero scalars. We will use the notation $\tau$ introduced in (\ref{datum_2_tau}).

Our goal now is to determine the parameters $\mu_i\in\FF^\times$ that appear in the above formula. On the one hand, the automorphism $\vp^2$ is the conjugation by matrix ${}^t\Phi^{-1}\Phi$ given by
\[
{}^t\Phi^{-1}\Phi=\sum_{i=1}^\ell\beta(t_i)I_{q_i}\ot I\oplus \sum_{i=\ell+1}^m\beta(t_i)\sgn(S_i)I_{2q_i}\ot I\oplus
\sum_{i=m+1}^k\begin{pmatrix}\mu_i I_{q_i}&0\\0&\mu_i^{-1} I_{q_i}\end{pmatrix}\ot I.
\]
On the other hand, $\vp^2$ acts as $\chi^2$. We now derive the conditions that are necessary and sufficient for $\chi^2\ast X=({}^t\Phi^{-1}\Phi{})^{-1}X({}^t\Phi^{-1}\Phi)$ to hold for all $X\in R$. 

Recall the idempotents $e_1,\ld,e_k\in C$ defined earlier (where $e_i=e_i'+e_i''$ for $i>m$). We denote by $U_{ij}$ any matrix in the Peirce component $e_iCe_j$. Then, for $1\le i,j\le m$, we have, for all $\bt\in\bT$,
\[
\chi^2\ast (U_{ij}\ot X_{\bt})=\chi^2(\bg_i^{-1}\bg_j\bt)U_{ij}\ot X_{\bt} 
\]
while 
\[
\begin{split}
\vp^2(U_{ij}\ot X_{\bt})&=({}^t\Phi^{-1}\Phi{})^{-1}(U_{ij}\ot X_{\bt})({}^t\Phi^{-1}\Phi)\\
&=\beta(\bt_i)\sgn(S_i)\beta(\bt_j)\sgn(S_j)(U_{ij}\ot X_{\bt}).
\end{split}
\]
For $m+1\le i,j\le k$, we write $U_{ij}=\matr{A&B\\C&D}$ according to the decompositions $e_i=e_i'+e_i''$ and $e_j=e_j'+e_j''$. Then, for all $\bt\in\bT$, 
\[
\chi^2\ast(U_{ij}\ot X_\bt)
=\begin{pmatrix}\chi^2((\bg_i')^{-1}\bg_j')A&\chi^2((\bg_i')^{-1}\bg_j'')B\\
\chi^2((\bg_i'')^{-1}\bg_j')C&\chi^2((\bg_i'')^{-1}\bg_j'')D\end{pmatrix},
\]
while
\[
\vp^2(U_{ij}\ot X_{\bt})=\begin{pmatrix}\mu_i^{-1}\mu_j A&\mu_i^{-1}\mu_j^{-1}B\\ 
\mu_i\mu_j C&\mu_i\mu_j^{-1} D\end{pmatrix}.
\]
For $1\le i\le m$ and $m+1\le j\le k$, we write $U_{ij}=\begin{pmatrix}A&B\end{pmatrix}$ according to the decomposition $e_j=e_j'+e_j''$. Then, for all $\bt\in\bT$,
\[
\chi^2\ast(U_{ij}\ot X_\bt)=\begin{pmatrix}\chi^2(\bg_i^{-1}\bg_j')A&\chi^2(\bg_i^{-1}\bg_j'')B\end{pmatrix},
\]
while
\[
\vp^2(U_{ij}\ot X_{\bt})=\beta(\bt_i)\sgn(S_i)\begin{pmatrix}\mu_j A&\mu_j^{-1}B\end{pmatrix}.
\]
For $m+1\le i\le k$ and $1\le j\le m$, we have a similar calculation.
 
By way of comparison, we derive $\chi^2(\bt)=\mathrm{const}$ for all $\bt\in\bT$, and so $\chi^2(\bT)=1$. Hence the natural epimorphism $\pi: G\to\bG$ splits over $\bT$, i.e., $\pi^{-1}(\bT)=T\times\langle h\rangle$, where $T=\pi^{-1}(\bT)\cap\ker\chi$. So we may identify $T$ with $\bT$ and write $t_i$ for the representative of the coset $\bt_i$ in $T$. Conversely, if $\pi: G\to\bG$ splits over $\bT$, then $\chi^2(\bT)=1$.

In the case $1\le i,j\le m$, our relations are equivalent to $\beta(t_i)\sgn(S_i)\chi^2(\bg_i)=\beta(t_j)\sgn(S_j)\chi^2(\bg_j)$. Therefore, we have a fixed $\lambda\in \FF^\times$ such that 
\bee{eq1}
\beta(t_i)\sgn(S_i)\chi^2(\bg_i)=\lambda\quad\mbox{for all}\quad i=1,\ld,m.
\ene

In the case $m+1\le i, j\le k$, our relations are equivalent to  
\[
\mu_i^{-1}\chi^2(\bg_i')=\mu_j^{-1}\chi^2(\bg_j')
\]
and
\[
\mu_i^{-1}\chi^2(\bg_i')=\mu_j\chi^2(\bg_j'').
\]
Therefore, we have a fixed $\mu\in\FF^\times$ such that 
\bee{eq2}
\mu_i^{-1}\chi^2(\bg_i')=\mu_i\chi^2(\bg_i'')=\mu\quad\mbox{for all}\quad i=m+1,\ld,k.
\ene

In the case $1\le i\le m$ and  $m+1\le j\le k$, our relations are equivalent to 
\bee{eq3}
\mu_j^{-1}\chi^2(\bg_j')=\beta(t_i)\sgn(S_i)\chi^2(\bg_i)=\mu_j\chi^2(\bg_j'').
\ene

If both (\ref{eq1}) and (\ref{eq2}) are present (i.e., $m\neq 0,k$), then (\ref{eq3}) is equivalent to $\mu=\lambda$. We have proved that if the $\bG$-grading on $R$ is the coarsening a $G$-grading on $\Lie{R}$ induced by $\pi:G\to\bG$, and $\chi$ acts on $R$ as $-\vp$, then $\pi^{-1}(\bT)$ splits and conditions (\ref{eq1}) and (\ref{eq2}) hold with $\lambda=\mu$. Conversely, if $R=\M(\bG,\bT,\beta,\kappa,\gamma)$ is such that $\pi^{-1}(\bT)$ splits, and an anti-automorphism $\vp$ is given by matrix (\ref{formula_antaut_restated}) such that (\ref{eq1}) and (\ref{eq2}) hold with $\lambda=\mu$, then $\vp^2$ acts as $\chi^2$ on $R$ and hence $-\vp$ defines a refinement of the $\bG$-grading on $R$ to a $G$-grading (as a vector space). The latter is automatically a grading of the Lie algebra $\Lie{R}$, since $-\vp$ is an automorphism of $\Lie{R}$.

To summarize, we state the following

\begin{prop}\label{chi_phi}
Let $h\in G$, $\pi:G\to\bG=G/\langle h \rangle$ and $R=\M(\bG,\bT,\beta,\kappa,\gamma)$ be as above. Let $H=\pi^{-1}(\bT)$. Fix $\chi\in\wh{G}$ with $\chi(h)=-1$. Let $\vp$ be the anti-automorphism of the $\bG$-graded algebra $R$ given by $\vp(X)=\Phi^{-1}({}^tX)\Phi$ with $\Phi$ as in (\ref{formula_antaut_restated}). Then
\[
R_g=\{X\in R_\bg\;|\;-\vp(X)=\chi(g)X\}\quad\mbox{for all}\quad g\in G
\]
defines a $G$-grading on $\Lie{R}$ if and only if $H$ splits as $T\times\langle h\rangle$ with $T=H\cap\ker\chi$ and the following condition holds (identifying $\bT$ with $T$):
\begin{equation}\label{compat_chi}
\begin{split}
&\beta(t_1)\chi^2(\bg_1)=\ldots=\beta(t_\ell)\chi^2(\bg_\ell)\\
&=\beta(t_{\ell+1})\sgn(S_{\ell+1})\chi^2(\bg_{\ell+1})=\ldots=\beta(t_m)\sgn(S_m)\chi^2(\bg_m)\\
&=\mu_{m+1}^{-1}\chi^2(\bg_{m+1}')=\mu_{m+1}\chi^2(\bg_{m+1}'')=\ldots=\mu_k^{-1}\chi^2(\bg_k')=\mu_k\chi^2(\bg_k'').
\end{split}
\end{equation}
\end{prop}\hfill{$\square$}

It is convenient to distinguish the following three cases for a grading of Type II on $L$:
\begin{itemize}
\item The case with $\ell>0$ will be referred to as Type $\mathrm{II}_1$;
\item The case with $\ell=0$ but $m>0$, will be referred to as Type $\mathrm{II}_2$;
\item The case with $m=0$ will be referred to as Type $\mathrm{II}_3$.
\end{itemize}

\begin{df}
We will say that $\gamma$ is {\em admissible} if it satisfies
\begin{equation}\label{compat_A_gamma_mod} 
\bg_1^2\equiv\ldots\equiv\bg_m^2\equiv\bg_{m+1}'\bg_{m+1}''\equiv\ldots\equiv\bg_k'\bg_k''\pmod{\bT}
\end{equation}
and, for some $\bt_1,\ldots,\bt_\ell\in \bT$, we have $\bg_1^2\bt_1=\ldots=\bg_\ell^2\bt_\ell$ and
\begin{equation}\label{admissible}
\beta(\bt_1)\chi^2(\bg_1)=\ldots=\beta(\bt_\ell)\chi^2(\bg_\ell).
\end{equation}
(If $\ell\leq 1$, then condition (\ref{admissible}) is automatically satisfied.) 
\end{df}

Note that the above definition does not depend on the choice of $\chi\in\wh{G}$ with $\chi(h)=-1$. Indeed, if we replace $\chi$ by $\wt{\chi}=\chi\psi$ where $\psi\in\wh{G}$ satisfies $\psi(h)=1$, then $\psi$ can be regarded as a character on $\bG$ and we can compute: 
\[
\wt{\chi}^2(\bg_i^{-1}\bg_j)=\chi^2(\bg_i^{-1}\bg_j)\psi^2(\bg_i^{-1}\bg_j)=\chi^2(\bg_i^{-1}\bg_j)\psi(\bg_i^{-2}\bg_j^2)=\chi^2(\bg_i^{-1}\bg_j)\psi(\bt_i\bt_j)
\]
for all $1\leq i,j\leq\ell$. On the other hand, for $\bt\in \bT$, we have 
\[
\beta(\bt\bt_i)\beta(\bt\bt_j)=\beta(\bt)\beta(\bt_i)\beta(\bt,\bt_i)\beta(\bt)\beta(\bt_j)\beta(\bt,\bt_j)
=\beta(\bt_i)\beta(\bt_j)\beta(\bt,\bt_i\bt_j).
\]
Therefore, if condition (\ref{admissible}) holds for $\chi$ and $\bt_1,\ldots,\bt_\ell$, then it holds for $\wt{\chi}$ and $\bt\bt_1,\ldots,\bt\bt_\ell$ where $\bt$ is the unique element of $\bT$ such that $\beta(\bt,\overline{u})=\psi(\overline{u})$ for all $\overline{u}\in\bT$.

As pointed out earlier, for $\gamma$ satisfying (\ref{compat_A_gamma_mod}), we can replace $\bg_i''$, $i>m$, within their cosets mod $\bT$ so that $\gamma$ satisfies (\ref{compat_A_gamma}). 

We now give our standard realizations for gradings of Type II. Let $H\subset G$ be an elementary $2$-group of odd rank containing $h$. Let $\beta$ be a nondegenerate alternating bicharacter on $\bT=H/\langle h\rangle$. Fix $\kappa$. Choose $\gamma$ formed from elements of $\bG=G/\langle h\rangle$ and $\tau$ formed from elements of $\bT=H/\langle h\rangle$ so that they satisfy (\ref{compat_A_gamma}). Let $R=\M(\bG,\bT,\beta,\kappa,\gamma)$. Fix $\chi\in\wh{G}$ with $\chi(h)=-1$ and identify $\bT$ with $T=H\cap\ker\chi$.

\begin{df}\label{type_II_1}
Suppose $\ell>0$ and $\gamma$ is admissible. Let $\Phi$ be the matrix given by (\ref{formula_antaut_restated}) where the scalars $\mu_i$ and matrices $S_i$ are determined by equation (\ref{compat_chi}). Then, by Proposition \ref{chi_phi}, the anti-automorphism $\vphi(X)=\Phi^{-1}({}^tX)\Phi$ defines a refinement of the $\bG$-grading on the associative algebra $R$ to a $G$-grading $R=\bigoplus_{g\in G}R_g$ as a Lie algebra. Set $L_g=(R_g\cap[R,R])$ mod $Z(R)$. We will denote the $G$-graded algebra $L$ obtained in this way as $\AII{1}(G,H,h,\beta,\kappa,\gamma,\tau)$.
\end{df}

\begin{df}\label{type_II_2}
Suppose $\ell=0$ and $m>0$. Choose $\delta=(\delta_1,\ldots,\delta_m)\in\{\pm 1\}^m$ so that
\[
\beta(t_1)\chi^2(\bg_1)\delta_1=\ldots=\beta(t_m)\chi^2(\bg_m)\delta_m.
\]
(Note that there are exactly two such choices.) Let $\Phi$ be the matrix given by (\ref{formula_antaut_restated}) where the matrices $S_i$ are selected by the rule $\sgn(S_i)=\delta_i$ and the scalars $\mu_i$ are determined by equation (\ref{compat_chi}). Then, by Proposition \ref{chi_phi}, the anti-automorphism $\vphi(X)=\Phi^{-1}({}^tX)\Phi$ defines a refinement of the $\bG$-grading on the associative algebra $R$ to a $G$-grading $R=\bigoplus_{g\in G}R_g$ as a Lie algebra. Set $L_g=(R_g\cap[R,R])$ mod $Z(R)$. We will denote the $G$-graded algebra $L$ obtained in this way as $\AII{2}(G,H,h,\beta,\kappa,\gamma,\tau,\delta)$.
\end{df}

\begin{df}\label{type_II_3}
Suppose $m=0$. Then we have
\[
\chi^2(\bg_1'\bg_1'')=\ldots=\chi^2(\bg_k'\bg_k'').
\]
Let $\mu$ be a scalar such that $\mu^2$ is equal to the common value of $\chi^2(\bg_i'\bg_i'')$. (There are two choices.) Let $\Phi$ be the matrix given by (\ref{formula_antaut_restated}) where the scalars $\mu_i$ are determined by equation 
\[
\mu_{1}^{-1}\chi^2(\bg_{1}')=\mu_{1}\chi^2(\bg_{1}'')=\ldots=\mu_k^{-1}\chi^2(\bg_k')=\mu_k\chi^2(\bg_k'')=\mu.
\]
Then, by Proposition \ref{chi_phi}, the anti-automorphism $\vphi(X)=\Phi^{-1}({}^tX)\Phi$ defines a refinement of the $\bG$-grading on the associative algebra $R$ to a $G$-grading $R=\bigoplus_{g\in G}R_g$ as a Lie algebra. Set $L_g=(R_g\cap[R,R])$ mod $Z(R)$. We will denote the $G$-graded algebra $L$ obtained in this way as $\AII{3}(G,H,h,\beta,\kappa,\gamma,\mu)$.
\end{df}

\begin{df}\label{equiv_II_2}
Referring to Definition \ref{type_II_2}, we will write $(\kappa,\gamma,\tau,\delta)\approx (\wt{\kappa},\wt{\gamma},\wt{\tau},\wt{\delta})$ if $\kappa$ and $\wt{\kappa}$ have the same number of components of each type, i.e., the same values of $m$ and $k$, and there exist an element $\bg\in \bG$ and a permutation $\pi$ of the symbols $\{1,\ldots,k\}$ preserving the sets $\{1,\ldots,m\}$ and $\{m+1,\ldots,k\}$ such that $\wt{q}_i=q_{\pi(i)}$ for all $i$, $\wt{t}_i=t_{\pi(i)}$, $\wt{\bg_i}\equiv \bg_{\pi(i)}\bg\pmod{\bT}$ and $\wt{\delta}_i=\delta_{\pi(i)}$ for all $i=1,\ldots,m$, and $\{\wt{\bg'_i},\wt{\bg''_i}\}\equiv\{\bg'_{\pi(i)}\bg,\bg''_{\pi(i)}\bg\}\pmod{\bT}$ for all $i=m+1,\ldots,k$.
\end{df}

\begin{df}\label{equiv_II_3}
Referring to Definition \ref{type_II_3}, we will write $(\kappa,\gamma,\mu)\approx (\wt{\kappa},\wt{\gamma},\wt{\mu})$ if $\kappa$ and $\wt{\kappa}$ have the same number of components $k$ and there exist an element $\bg\in\bG$ and a permutation $\pi$ of the symbols $\{1,\ldots,k\}$ such that $\wt{q}_i=q_{\pi(i)}$, $\{\wt{\bg'_i},\wt{\bg''_i}\}\equiv\{\bg'_{\pi(i)}\bg,\bg''_{\pi(i)}\bg\}\pmod{\bT}$ and  $\wt{\bg'_i}\wt{\bg''_i}=\bg'_{\pi(i)}\bg''_{\pi(i)}\bg^2$ for all $i$, and, finally, $\wt{\mu}=\mu\chi^2(\bg)$.
\end{df}

\begin{theorem}\label{A_main}
Let $\FF$ be an algebraically closed field, $\chr{\FF}\neq 2$. Let $G$ be an abelian group. Let $L=\Psl_n(\FF)$ where $n\geq 3$. If $n=3$, assume also that $\chr{\FF}\neq 3$. Let $L=\bigoplus_{g\in G}L_g$ be a $G$-grading. Then the graded algebra $L$ is isomorphic to one of the following:
\begin{itemize}
\item $\AI(G,T,\beta,\kappa,\gamma)$,
\item $\AII{1}(G,H,h,\beta,\kappa,\gamma,\tau)$,
\item $\AII{2}(G,H,h,\beta,\kappa,\gamma,\tau,\delta)$,
\item $\AII{3}(G,H,h,\beta,\kappa,\gamma,\mu)$,
\end{itemize}
as in Definitions \ref{type_I}, \ref{type_II_1}, \ref{type_II_2} and \ref{type_II_3}, with $|\kappa|\sqrt{|T|}=n$ in Type I and $|\kappa|\sqrt{|H|/2}=n$ in Type II.
Graded algebras belonging to different types listed above are not isomorphic. Within each type, we have the following:
\begin{itemize}
\item $\AI(G,T_1,\beta_1,\kappa_1,\gamma_1)\cong\AI(G,T_2,\beta_2,\kappa_2,\gamma_2)$ if and only if $T_1=T_2$, $\beta_1=\beta_2$, and  $(\kappa_1,\gamma_1)\sim(\kappa_2,\gamma_2)$ or $(\kappa_1,\gamma_1)\sim(\kappa_2,\gamma_2^{-1})$;
\item $\AII{1}(G,H_1,h_1,\beta_1,\kappa_1,\gamma_1,\tau_1)\cong\AII{1}(G,H_2,h_2,\beta_2,\kappa_2,\gamma_2,\tau_2)$ if and only if $H_1=H_2$, $h_1=h_2$, $\beta_1=\beta_2$, and $(\kappa_1,\gamma_1,\tau_1)\approx(\kappa_2,\gamma_2,\tau_2)$ or $(\kappa_1,\gamma_1,\tau_1)\approx(\kappa_2,\gamma_2^{-1},\tau_2)$;
\item $\AII{2}(G,H_1,h_1,\beta_1,\kappa_1,\gamma_1,\tau_1,\delta_1)\cong\AII{2}(G,H_2,h_2,\beta_2,\kappa_2,\gamma_2,\tau_2,\delta_2)$ if and only if $H_1=H_2$, $h_1=h_2$, $\beta_1=\beta_2$, and $(\kappa_1,\gamma_1,\tau_1,\delta_1)\approx(\kappa_2,\gamma_2,\tau_2,\delta_2)$ or $(\kappa_1,\gamma_1,\tau_1,\delta_1)\approx(\kappa_2,\gamma_2^{-1},\tau_2,\delta_2)$;
\item $\AII{3}(G,H_1,h_1,\beta_1,\kappa_1,\gamma_1, \mu_1)\cong\AII{3}(G,H_2,h_2,\beta_2,\kappa_2,\gamma_2,\mu_2)$ if and only if $H_1=H_2$, $h_1=h_2$, $\beta_1=\beta_2$,  and $(\kappa_1,\gamma_1,\mu_1)\approx(\kappa_2,\gamma_2,\mu_2)$ or $(\kappa_1,\gamma_1,\mu_1)\approx(\kappa_2,\gamma_2^{-1},\mu_2^{-1})$.
\end{itemize}
\end{theorem}

\begin{proof}
The first statement is a combination of Theorem \ref{form_antaut} and Proposition \ref{chi_phi}. The non-isomorphism of graded algebras belonging to different types is clear.

For Type I, let $R_1=\M(G,T_1,\beta_1,\kappa_1,\gamma_1)$ and $R_2=(G,T_2,\beta_2,\kappa_2,\gamma_2)$. By Theorem \ref{matr_main}, $R_1\cong R_2$ if and only if $T_1=T_2$, $\beta_1=\beta_2$, and  $(\kappa_1,\gamma_1)\sim(\kappa_2,\gamma_2)$. It remains to observe that the outer automorphism $X\mapsto -{}^t X$ transforms $\M(G,T,\beta,\kappa,\gamma)$ to $\M(G,T,\beta,\kappa,\gamma^{-1})$.

For Type II, the element $h$, the subgroup $H$, and the bicharacter $\beta$ on $\bT=H/\langle h\rangle$ are uniquely determined by the grading, so we may assume $H_1=H_2$,  $h_1=h_2$, and $\beta_1=\beta_2$. Let $R_1=\M(\bG,\bT,\beta,\kappa_1,\gamma_1)$ and $R_2=\M(\bG,\bT,\beta,\kappa_2,\gamma_2)$. Fix $\chi\in\wh{G}$ with $\chi(h)=-1$.  Let $\vp_1$ and $\vp_2$ be the corresponding anti-automorphisms. We have to check that $(R_1,\vp_1)\cong(R_2,\vp_2)$ if and only if 
\begin{enumerate}
\item[$\mathrm{II}_1$)] $(\kappa_1,\gamma_1,\tau_1)\approx(\kappa_2,\gamma_2,\tau_2)$,
\item[$\mathrm{II}_2$)] $(\kappa_1,\gamma_1,\tau_1,\delta_1)\approx(\kappa_2,\gamma_2,\tau_2,\delta_2)$,
\item[$\mathrm{II}_3$)] $(\kappa_1,\gamma_1,\mu_1)\approx(\kappa_2,\gamma_2,\mu_2)$.
\end{enumerate}

For Type $\mathrm{II}_1$, the ``only if'' part is clear, since $(\kappa,\gamma,\tau)$ is an invariant of $(R,\vp)$ (up to transformations indicated in the definition of the equivalence relation $\approx$). Indeed, $(\kappa,\gamma)$ is an invariant of the $\bG$-grading, and $\tau$ corresponds to the restrictions of $\vp$ to $D_1,\ldots,D_m$ by Proposition \ref{antaut_division}. To prove the ``if'' part, assume $(\kappa_1,\gamma_1,\tau_1)\approx(\kappa_2,\gamma_2,\tau_2)$. Then, by Theorem \ref{matr_main}, there exists an isomorphism of $\bG$-graded algebras $\psi:R_1\to R_2$. By Remark \ref{monomial_iso}, we can take for $\psi$ a monomial isomorphism associated to the permutation $\pi$ in Definition \ref{equiv_2}. The matrix of the anti-automorphism $\psi^{-1}\vphi_2\psi$ on $R_1$ is then obtained from the matrix of $\vphi_2$ by permuting the blocks on the diagonal so that they align with the corresponding blocks of $\vphi_1$, and possibly multiplying some of the blocks by $-1$. Hence, by Theorem \ref{form_antaut}, $\psi^{-1}\vphi_2\psi$ can be transformed to $\vphi_1$ by an automorphism of the $\bG$-graded algebra $R_1$.

For Type $\mathrm{II}_2$, the proof is similar, since $\delta$ corresponds to the restrictions of $\vp$ to $C_1,\ldots,C_m$ and thus is an invariant of $(R,\vp)$.

For Type $\mathrm{II}_3$, we show in the same manner that if $(\kappa_1,\gamma_1,\mu_1)\approx(\kappa_2,\gamma_2,\mu_2)$, then $(R_1,\vp_1)\cong(R_2,\vp_2)$. Namely, we take a monomial isomorphism of $\bG$-graded algebras $\psi:R_1\to R_2$ associated to the permutation $\pi$ in Definition \ref{equiv_II_3}. The effect of $\psi$ on  $\Phi_2$ is just the permutation of blocks. 
The factor $\chi^2(\bg)$ in Definition \ref{equiv_II_3} makes sure that the block with $\mu_i=\mu^{-1}\chi^2(\bg'_i)$ in $\Phi_2$ matches up with the block with $\mu_{\pi(i)}=\mu^{-1}\chi^2(\bg'_{\pi(i)})$ in $\Phi_1$. Conversely, suppose there exists an isomorphism $\psi:(R_1,\vp_1)\to(R_2,\vp_2)$. As in the proof of Corollary \ref{matr_with_inv}, we write $\psi=\psi_0\alpha$ where $\psi_0$ is a monomial isomorphism and $\alpha$ is in $\PGL_{\kappa_1}(\FF)\times\Aut_G(D)$. The action of $\psi_0$ on $\vp_2$ permutes the blocks and replaces the second tensor factor $I$ by $X_{t_0}$ for some $t_0\in T$. The action of $\alpha$ on $\psi_0^{-1}\vp_2\psi_0$ cannot change $t_0$ or the values of the scalars. 
We conclude that $t_0=e$ and $(\kappa_1,\gamma_1,\mu_1)\approx(\kappa_2,\gamma_2,\mu_2)$.
\end{proof}



\begin{remark}
Let $\FF$ and $G$ be as in Theorem \ref{A_main}. Let $L=\Sl_2(\FF)$. If $L=\bigoplus_{g\in G}L_g$ is a $G$-grading, then the graded algebra $L$ is isomorphic to $\AI(G,T,\beta,\kappa,\gamma)$ where $|\kappa|\sqrt{|T|}=2$. This, of course, gives two possibilities: either $T=\{e\}$ or $T\cong\ZZ_2^2$. In the first case the $G$-grading is induced from a Cartan decomposition by a homomorphism $\ZZ\to G$. The isomorphism classes of such gradings are in one-to-one correspondence with unordered pairs of the form  $\{g,g^{-1}\}$, $g\in G$. In the second case the $G$-grading is given by Pauli matrices. The isomorphism classes of such gradings are in one-to-one correspondence with subgroups $T\subset G$ such that $T\cong\ZZ_2^2$.
\end{remark}

\begin{remark}
The remaining case $L=\Psl_3(\FF)$ where $\chr{\FF}=3$ can be handled using octonions. Let $\OO$ be the algebra of octonions over an algebraically closed field $\FF$. Then the subspace $\OO'$ of zero trace octonions is a Malcev algebra with respect to the commutator $[x,y]=xy-yx$. If $\chr{\FF}=3$, then $\OO'$ is a Lie algebra isomorphic to $L$. Assuming $\chr{\FF}\neq 2$, we have $xy=\frac{1}{2}([x,y]-n(x,y)1)$ for all $x,y\in\OO'$, where $n$ is the norm of $\OO$. We also have $(\ad x)^3=-4n(x)(\ad x)$ for all $x\in\OO'$. It follows that if $\psi$ is an automorphism of $\OO'$, then $\psi$ preserves $n$ and, setting $\psi(1)=1$, we obtain an automorphism of $\OO$. Hence the restriction map $\Aut(\OO)\to\Aut(\OO')$ is an isomorphism of algebraic groups. Similarly, one shows that the restriction map $\Der(\OO)\to\Der(\OO')$ is an isomorphism of Lie algebras.\footnote{This argument was communicated to us by A. Elduque.} It follows that $\AAut(\OO')$ is smooth and can be identified with the algebraic group $\Aut(\OO)$. In particular, this means that the isomorphism classes of $G$-gradings on $\OO$ are in one-to-one correspondence (via restriction) with the isomorphism classes of $G$-gradings on $\OO'$ (cf. \cite[Theorem 9]{Eld98}). 

All gradings on $\OO$ (in any characteristic) were described in \cite{Eld98}. For $\chr{\FF}\neq 2$, they are of two types: 
\begin{itemize}
\item ``elementary'' gradings obtained by choosing $g_1,g_2,g_3\in G$ with $g_1g_2g_3=e$ and assigning degree $e$ to $e_1$ and $e_2$, degree $g_i$ to $u_i$ and degree $g_i^{-1}$ to $v_i$, $i=1,2,3$, where $\{e_1,e_2,u_1,u_2,u_3,v_1,v_2,v_3\}$ is a canonical basis for $\OO$;
\item ``division'' gradings by $\ZZ_2^3$ obtained by iterating the Cayley-Dickson doubling process three times.
\end{itemize}
It is easy to see when two $G$-gradings on $\OO$ are isomorphic. The isomorphism classes of ``elementary'' gradings are in one-to-one correspondence with unordered pairs of the form  $\{S,S^{-1}\}$ where $S$ is an unordered triple $\{g_1,g_2,g_3\}$, $g_i\in G$ with $g_1g_2g_3=e$. The isomorphism classes of ``division'' gradings are in one-to-one correspondence with subgroups $T\subset G$ such that $T\cong\ZZ_2^3$. An ``elementary'' grading is not isomorphic to a ``division'' grading.

If $\chr{\FF}=3$, then the above is also the classification of $G$-gradings on $L=\Psl_3(\FF)$. As shown in \cite{Ksur}, up to isomorphism, any grading on $L$ is induced from the matrix algebra $M_3(\FF)$. Namely, any ``elementary'' grading on $L$ can be obtained as a Type I grading, and any ``division'' grading on $L$ is isomorphic to a Type II gradings. The only difference with the case of $\Sl_3(\FF)$ where $\chr{\FF}\neq 3$ is that there are fewer isomorphism classes of gradings in characteristic 3 (in particular, some ``Type II'' gradings are isomorphic to ``Type I'' gradings).
\end{remark}


\section{Gradings on Lie Algebras of Types $\B$, $\C$, $\D$}\label{section_BCD}

The classification of gradings for Lie algebras $\So_n(\FF)$ and $\Sp_n(\FF)$ follows immediately from Corollary \ref{matr_with_inv}. We state the results here for completeness. Recall $\M^*(G,T,\beta,\kappa,\gamma,\tau,\delta)$ from Definition \ref{datum_2}. Let $L=\sks(R,\vphi)=\{X\in R\;|\;\vphi(X)=-X\}$. Then $L=\bigoplus_{g\in G}L_g$ where $L_g=R_g\cap L$.

\begin{df}\label{BCD} 
Let $n=|\kappa|\sqrt{|T|}$.
\begin{itemize}
\item If $\delta=1$ and $n$ is odd, then necessarily $T=\{e\}$. We will denote the $G$-graded algebra $L$ by $\B(G,\kappa,\gamma)$.
\item If $\delta=-1$ (hence $n$ is even), then we will denote the $G$-graded algebra $L$ by $\C(G,T,\beta,\kappa,\gamma,\tau)$.
\item If $\delta=1$ and $n$ is even, then we will denote the $G$-graded algebra $L$ by $\D(G,T,\beta,\kappa,\gamma,\tau)$.
\end{itemize}
\end{df}

\begin{theorem}\label{BCD_main}
Let $\FF$ be an algebraically closed field, $\chr{\FF}\neq 2$. Let $G$ be an abelian group. 
\begin{itemize}
\item Let $L=\So_n(\FF)$, with odd $n\geq 5$. Let $L=\bigoplus_{g\in G}L_g$ be a $G$-grading. Then the graded algebra $L$ is isomorphic to $\B(G,\kappa,\gamma)$,
\item Let $L=\Sp_n(\FF)$, with even $n\geq 6$. Let $L=\bigoplus_{g\in G}L_g$ be a $G$-grading. Then the graded algebra $L$ is isomorphic to $\C(G,T,\beta,\kappa,\gamma,\tau)$,
\item Let $L=\So_n(\FF)$, with even $n\geq 10$. Let $L=\bigoplus_{g\in G}L_g$ be a $G$-grading. Then the graded algebra $L$ is isomorphic to $\D(G,T,\beta,\kappa,\gamma,\tau)$, 
\end{itemize}
as in Definition \ref{BCD}. Also, under the above restrictions on $n$, we have the following:
\begin{itemize}
\item $\B(G,\kappa_1,\gamma_1)\cong\B(G,\kappa_1,\gamma_1)$ if and only if $(\kappa_1,\gamma_1)\approx(\kappa_1,\gamma_1)$;
\item $\C(G,T_1,\beta_1,\kappa_1,\gamma_1,\tau_1)\cong\C(G,T_2,\beta_2,\kappa_2,\gamma_2,\tau_2)$ if and only if $T_1=T_2$, $\beta_1=\beta_2$ and 
$(\kappa_1,\gamma_1,\tau_1)\approx(\kappa_2,\gamma_2,\tau_2)$;
\item $\D(G,T_1,\beta_1,\kappa_1,\gamma_1,\tau_1)\cong\D(G,T_2,\beta_2,\kappa_2,\gamma_2,\tau_2)$ if and only if $T_1=T_2$, $\beta_1=\beta_2$ and 
$(\kappa_1,\gamma_1,\tau_1)\approx(\kappa_2,\gamma_2,\tau_2)$.
\end{itemize}
\end{theorem}
\hfill{$\square$}


\end{document}